\pdfoutput=1
\documentclass[twocolumn]{autart}

\usepackage{harvard}

\usepackage{graphicx}
\usepackage{subcaption}
\usepackage{amssymb,amsmath,mathrsfs}

\usepackage{color}

\newcommand{\real}{\mathbb{R}}

\newcommand{\nnum}{\mathbb{N}}

\newcommand{\LInfNorm}[1]{\left\|{#1}\right\|_{\mathcal{L}_{\infty}}}
\newcommand{\norm}[1]{\|{#1}\|}

\begin{document}

\begin{frontmatter}
\runtitle{Exponential Stability of Differential Repetitive Processes}  
\title{Exponential Stability of Nonlinear Differential Repetitive Processes with Applications to Iterative Learning Control\thanksref{footnoteinfo}}

\thanks[footnoteinfo]{This paper was not presented at any IFAC 
meeting. Corresponding author B.~Alt{\i}n. Tel. +1 (831) 459-2939.}

\author[Altin]{Berk Alt{\i}n}\ead{berkaltin@ucsc.edu},
\author[Barton]{Kira Barton}\ead{bartonkl@umich.edu}

\address[Altin]{Department of Computer Engineering, University of California, Santa Cruz, CA 95064}  
\address[Barton]{Department of Mechanical Engineering, University of Michigan, Ann Arbor, MI 48109}             


\begin{abstract}

This paper studies exponential stability properties of a class of two-dimensional~(2D) systems called differential repetitive processes~(DRPs). Since a distinguishing feature of DRPs is that the problem domain is bounded in the ``time'' direction, the notion of stability to be evaluated does not require the nonlinear system defining a DRP to be stable in the typical sense. In particular, we study a notion of exponential stability along the discrete iteration dimension of the 2D dynamics, which requires the boundary data for the differential pass dynamics to converge to zero as the iterations evolve. Our main contribution is to show, under standard regularity assumptions, that exponential stability of a DRP is equivalent to that of its linearized dynamics. In turn, exponential stability of this linearization can be readily verified by a spectral radius condition. The application of this result to Picard iterations and iterative learning control~(ILC) is discussed. Theoretical findings are supported by a numerical simulation of an ILC algorithm.

\end{abstract}

\end{frontmatter}


\section{Introduction}
\label{sec:introduction}
For recursive nonlinear systems in the explicit form
\begin{equation}
	\left\{
	\begin{aligned}
		\dot{x}_{k+1}(t)		&= f(x_{k+1}(t),y_{k}(t),t),\\
				y_{k+1}(t)		&= g(x_{k+1}(t),y_{k}(t),t),
	\end{aligned}
	\right.
\label{eq:nonlinear}
\end{equation}
where~$(t,k)\in [0,T]\times\{0,1,\dots\}$ for some~${T\in[0,\infty)}$, we are interested in finding necessary and sufficient conditions that establish local exponential stability. The vectors~${x_k(t)\in\mathbb{R}^n}$ and~${y_k(t)\in\mathbb{R}^m}$ of this model represent the state and output, respectively. To uniquely determine the solution of~\eqref{eq:nonlinear}, it will be necessary to specify boundary conditions~$y_0$ and~${\mathbf{x}(0)\triangleq\{x_{k+1}(0)\}_{k=0}^{\infty}}$.

Roughly speaking, the notions of stability to be studied throughout this paper will be weak, in the sense that they will not require the one-dimensional~(1D) control system given by~$f$ to be stable. For example, exponential stability of~\eqref{eq:nonlinear} will imply that the function sequence~$\{y_k\}_{k=0}^{\infty}$ converges exponentially to zero in an appropriate signal norm, provided the boundaries are small, and~$\mathbf{x}(0)$ also converges exponentially to zero. The precise meaning of stability for this class of systems will be defined later in Section~\ref{sec:state}.

The nonlinear system~\eqref{eq:nonlinear} appears in many practical problems of interest and falls into the larger class of two-dimensional~(2D) dynamic systems called repetitive~(or multipass, earlier in the literature) processes,\footnote{Not to be confused with repetitive control.} in which information propagation occurs along two axes of independent variables. These processes are characterized by a sequence of passes with \textit{finite length} that act as forcing functions on the dynamics of future passes~\cite{rogers}: The output solution sequence~${\{y_{k}\}_{k=0}^{\infty}}$ of~\eqref{eq:nonlinear} can be found by applying the nonlinear system with differential dynamics described by the functions~$f$ and~$g$ in a repetitive manner. Hence, we will call any system of the form~\eqref{eq:nonlinear} a \textit{differential repetitive process~(DRP)}. The counterpart of the DRP~\eqref{eq:nonlinear} in the broader 2D systems theory, where it is assumed that~$T=\infty$, will be called a 2D mixed continuous-discrete time system.

The repetitive process paradigm arises in the modeling of certain engineering applications such as long wall coal cutting~\cite{edwards} and metal rolling~\cite{foda,edwards3}. A rich set of examples to these systems can also be found on a more abstract level since recursive algorithms for 1D dynamic systems can be treated as repetitive processes; e.g. iterative solutions to nonlinear optimal control problems~\cite{zidek,gupta}, nonlinear inversion methods~\cite{devasia}, iterative estimation and control design~\cite{albertos}, or the constructive proof of the Picard-Lindel\"{o}f theorem. A well-known class of algorithms that can be expressed in the repetitive process framework is iterative learning control (ILC)~\cite{kurek,hladowski,ahn}, wherein the inverse image of a desired output under a 1D input-output system is constructed through a recurrence relation inducing pass-to-pass dynamics. This problem will be tackled in Section~\ref{sec:picardILC}.

The study of DRPs and other 2D systems bearing similarities with~\eqref{eq:nonlinear} has a long history, beginning with the Roesser and Fornasini-Marchesini models introduced in the 1970s~\cite{roesser,fornasini1,fornasini2}. In particular, stability and performance properties of DRPs and 2D mixed continuous-discrete time systems, along with corresponding control strategies, have been researched extensively, \textit{predominantly for linear time-invariant~(LTI) systems}--see~\cite{rogers,chesi1,chesi2} and references therein. On the other hand, the need to develop rigorous stability tests in the nonlinear systems context has been highlighted only very recently. Among these works,~\cite{yeganefar} present forward and converse Lyapunov theorems for nonlinear Roesser models, with extensions to the stochastic case given in~\cite{pakshin}, and a~2D Lyapunov function approach is employed to prove exponential stability of DRPs in~\cite{emelianov}. It is also worth noting that the DRP~\eqref{eq:nonlinear} can be viewed as an infinite-dimensional hybrid system~\cite{liu,barreiro,sun} by concatenating the passes; e.g. by letting~${x(\tau,k+1)\triangleq x_{k+1}(t)}$ with~${\tau=t+kT}$, subject to the periodic reset~${x(kT,k+1)=x_{k+1}(0)}$, where~$T$ plays the role of an inherent delay,~$\tau$ the ordinary time, and~$k$ the jump time/index. As this reset function would change based on the prespecified boundary condition~$\mathbf{x}(0)$ and lacks any other structure, we will not follow a hybrid systems approach in the ensuing analysis. See also~\cite{rogers} for DRP modeling of a class of delay differential equations.

\begin{figure}[tbhp]
	\centering
	\includegraphics[width=1\columnwidth]{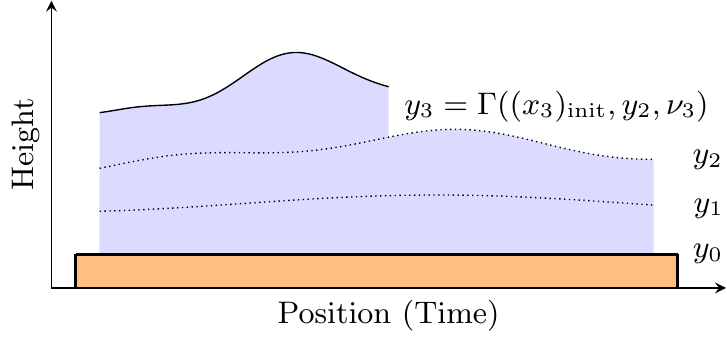}
	\caption{AM systems as repetitive processes: The substrate topography determines the initial output~$y_0$. The operator~$\Gamma$ maps the initial state~$(x_3)_{\mathrm{init}}$ and input~$\nu_3$ of pass~$3$~(in-layer dynamics), along with the prior pass profile~$y_2$~(layer-to-layer dynamics), to pass profile~$y_3$. The layer-to-layer dynamics is affected by physical phenomena such as material curing.}
	\label{fig:additive}
\end{figure}

The objective of this paper is to contribute to the recent literature on nonlinear repetitive process and 2D systems literature, and provide a connection between nonlinear DRPs of the form~\eqref{eq:nonlinear} and their linear counterparts. Therefore, our aim is to certify local exponential stability of DRPs via an appropriate linearization of~\eqref{eq:nonlinear}, and establish an analogue of the classical result that exponential stability of a 1D system is equivalent to that of its linear approximation, thereby expanding on the findings of~\cite{cdc2015}. Our primary motivation for this study comes from additive manufacturing~(AM) systems, wherein material in the fluid phase is often deposited in a layer-by-layer fashion~(Fig.~\ref{fig:additive}), leading to 2D dynamics: For instance, the laser metal deposition~(LMD) process is characterized by 1D (in-layer) dynamics that are height dependent due to heat transfer from prior layers~\cite{sammons}. It is possible to achieve accurate material distribution for the LMD process via linear repetitive process control techniques and a more control-oriented model consisting of static nonlinearities. This, however, requires the implicit assumption that the controlled nonlinear process is \textit{locally stable around its linearized equilibrium}~\cite{sammons3}. As a secondary motivation, in the ILC literature, it has been noted that nonlinear update laws have not been researched, save for adaptive laws for locally Lipschitz plants, and a systematic theory of nonlinear ILC is an open question~\cite{xusurvey,mooretutorial}.

The rest of the paper is organized as follows: Section~\ref{sec:state} introduces the necessary background, establishes the key Lipschitz property of the nonlinear operator, and states formal stability definitions.  Stability theory for LTI systems is extended to the linear time-varying~(LTV) case in Section~\ref{sec:linear}. Our main result, which establishes equivalence in terms of exponential stability between a DRP and its linearization, is presented in Section~\ref{sec:linearized}. Applications of this result to Picard iterations and ILC are discussed in Section~\ref{sec:picardILC}. An illustrative example is given in Section~\ref{sec:example} through an ILC system. Concluding remarks are given in Section~\ref{sec:conclusion}. In the hope of improving readability of the paper, proofs of certain technical results are given in Appendices~\ref{app:proofs},~\ref{app:claimproof} and~\ref{app:sufficiency}.

\section{Background and Preliminaries}
\label{sec:state}
This section will introduce the background material pertinent to our analysis, and lay out stability definitions for the DRP~\eqref{eq:nonlinear}. The precise definitions of stability to be presented will show the crucial difference between DRPs and 2D mixed continuous-discrete time systems, as the latter studies the trajectory of the real vector~${y_{k}(t)}$ over~${\{0,1,\dots\}\times[0,\infty)}$. In linear repetitive process theory, the gap between these two classes of systems is bridged via the stronger notion of \textit{stability along the pass}~\cite{rogers}, which requires the stability parameters to be~$T$ independent. Although this property is desirable in experimental implementations or numerical simulations, we will forgo this requirement for theoretical purposes.

\textbf{Notation}: We use~$\real$ to represent real numbers,~$\nnum$ nonnegative integers, and~$\mathbb{C}$ complex numbers. The spectral radius of a linear operator is denoted by~$\rho(.)$. The identity and zero operators are denoted as~$I$ and~0, respectively. For a real vector,~$\norm{.}_2$ is the~2 norm; in the rest of the paper~$\norm{.}$ will denote any of the equivalent norms in~$\real^p$.~$\mathcal{L}_p$ is the space of Lebesgue measurable functions on the compact interval~${[0,T]}$ with finite~$\mathcal{L}_p$ norm,~${p\in[1,\infty]}$. The space of all sequences on~$\real^p$ which converge to~0 is denoted as~$c_0$.

The inequalities below, stated without proof, will be of use for convergence analysis. {\color{black}Note that the convergence parameters~${2/(1-a)\geq 1}$ and~${(1+a)/2\in(0,1)}$ are continuous increasing functions of~$a$ on~$(0,1)$.}
\begin{claim}\label{claim:asymp}
Let~${\mathbf{a}\triangleq\{a_{k+1}\}_{k=0}^{\infty}}$ and~${\mathbf{b}\triangleq\{b_{k+1}\}_{k=1}^{\infty}}$ be real nonnegative sequences, where~$\mathbf{b}$ is bounded. Suppose that~${a_{k+1}=r a_k +b_{k+1}}$ for some~${r\in(0,1)}$ for all~${k\in\nnum}$. Then,~${\limsup_{k\to\infty}a_k\leq(1/(1-r))\limsup_{k\to\infty}b_k}$, and therefore~${\mathbf{b}\in c_0}$ implies~${\mathbf{a}\in c_0}$.
\end{claim}
\begin{claim}\label{claim:exp}
	Let~${a\in(0,1)}$. Then the sequence~${\{ka^{k-1}\}_{k=0}^{\infty}}$ is exponentially convergent and
	\[
		ka^{k-1}\leq\frac{2}{1-a}\left(\frac{1+a}{2}\right)^k, \quad \forall k\in\nnum.
	\]
\end{claim}


\subsection{The Nonlinear Operator over the Finite Horizon}

Before proceeding with further analysis, we will look at the properties of the system~\eqref{eq:nonlinear} as an input-state and input-output operator over the time interval~${[0,T]}$: Interchanging~$y_k$ with~$u$,~$x_{k+1}$ with~$\chi$, and~$y_{k+1}$ with~$w$, we consider 
\begin{equation}
	\left\{
	\begin{aligned}
		\dot{\chi}(t)		&= f(\chi(t),u(t),t),\\
				w(t)		&= g(\chi(t),u(t),t),
	\end{aligned}
	\right.
\label{eq:nonlinear1D}
\end{equation}
for all~${t\in[0,T]}$. The input~$u$ resides in~$\mathcal{Y}$, the space of continuously differentiable functions on~${[0,T]}$. We will impose the following standing assumptions on the nonlinear operator~$\Gamma$ that maps the pair~${(\chi(0),u)}$ to~$\chi$ and~$w$:
\begin{assum}
\label{assum:nonlinear}
	The nonlinear system~\eqref{eq:nonlinear1D} satisfies the following conditions:
	\begin{enumerate}
		\item The functions~$f$ and~$g$ vanish at the origin uniformly in time. That is,~${f(0,0,t)=0}$ and~${g(0,0,t)=0}$ for all~${t\in [0,T]}$.
		\item There exists~${\delta>0}$ such that for every~${(\chi(0),u)}$ that satisfies~${\norm{\chi(0)}+\LInfNorm{u}<\delta}$, there is a unique integral curve~$\chi$ of~\eqref{eq:nonlinear1D}, and~$\chi(t)$ is contained in a bounded open connected set~$X$ for all~${t\in[0,T]}$.
		\item There exists a compact set~${Y\subset\mathbb{R}^m}$ that contains the origin in its interior such that~$f$ and~$g$ are continuously differentiable in~${Z\triangleq \mathrm{cl}(X)\times Y\times [0,T]}$, where~$\mathrm{cl}(X)$ is the closure of~$X$.
	\end{enumerate}
\end{assum}

Assumption~\ref{assum:nonlinear} is a mild constraint on the system that bypasses the stability requirement in the time domain. We note that since~0 is an equilibrium of the differential equation, the set~$X$ must contain the origin. Without loss of generality, we will also assume that~$\delta$ is small enough so that~${\chi(0)\in X}$ and~${u(t)\in Y}$ for all~${t\in[0,T]}$ when~${\norm{\chi(0)}+\LInfNorm{u}<\delta}$. We denote by~$\Gamma_x$ the mapping~${(u,\chi(0))\mapsto \chi}$, and by~$\Gamma_y$ the mapping~${(u,\chi)\mapsto w}$, so
\[
	(w,\chi)=\Gamma(u,\chi(0))\triangleq( \Gamma_y(u,\Gamma_x(u,\chi(0))),\Gamma_x(u,\chi(0))).
\]
Now, we can show Lipschitz continuity of the operator~$\Gamma$ in the uniform norm topology. See Appendix~\ref{app:proofs} for a proof of this result.
\begin{lem}
	\label{lem:lipschitz}
	The nonlinear operator~$\Gamma$ given by~\eqref{eq:nonlinear1D} is locally Lipschitz with respect to~$(\chi(0),u)$. That is, there exist positive constants~$\bar{\delta}$ and~$L$ such that if
	\[
		(w_i,\chi_i)\triangleq\Gamma(\chi_i(0),u_i) \text{ and } \norm{\chi_i(0)}+\LInfNorm{u_i}<\bar{\delta},
	\]
	where~${i\in\{1,2\}}$, then
	\begin{equation*}
		\begin{aligned}
			\LInfNorm{\chi_1-\chi_2}		&\leq L(\LInfNorm{u_1-u_2}+\norm{\chi_1(0)-\chi_2(0)}),\\
			\LInfNorm{w_1-w_2}					&\leq L(\LInfNorm{u_1-u_2}+\norm{\chi_1(0)-\chi_2(0)}).
		\end{aligned}
	\end{equation*}
\end{lem}


\subsection{Boundary Dependent Stability Definitions}
We will now lay out definitions of stability for DRPs. First, we need the following norm to characterize exponential initial state sequences for exponential stability, which is similar to the conventional time-weighted norm used in the ILC literature:
\begin{defn}
	Let~${\mathbf{b}\triangleq\{b_{k+1}\}_{k=0}^{\infty}}$ be a sequence on~$\real^p$. For any~$\lambda\in(0,1]$, the exponential~$\lambda$~(${{e}_{\lambda}}$) norm of~$\mathbf{b}$ is defined as ${\norm{\mathbf{b}}_{{e}_{\lambda}}\triangleq\sup_{k\in\nnum}\lambda^{-k}\norm{b_{k+1}}}$.
\end{defn}

We leave it to the reader to verify that~${{e}_{\lambda}}$, the vector space of all sequences on~$\real^p$ with finite~${{e}_{\lambda}}$ norm, i.e. the space of sequences on~$\real^p$ that converge geometrically to~0 with rate faster than or equal to~${\lambda}$, satisfies~${e_{\lambda}\subset c_0\subset e_1\equiv \ell_{\infty}}$, for all~${\lambda\in(0,1)}$. The~$e_{\lambda}$ norm also satisfies 1) \textit{the shift property},~$\norm{\mathbf{b}_{\kappa}}_{{e}_{\lambda}}\leq\lambda^\kappa\norm{\mathbf{b}}_{{e}_{\lambda}}$, where~${\mathbf{b}_{\kappa}\triangleq\{b_{k+1}\}_{k=\kappa}^{\infty}}$, given any~${\kappa\in\nnum}$, and 2) \textit{the~$\lambda$ property}~${\norm{.}_{e_{\lambda_2}}\leq \norm{.}_{e_{\lambda_1}}}$ when~${0<\lambda_1\leq\lambda_2\leq 1}$.
\begin{defn}
	The origin of the DRP~\eqref{eq:nonlinear} is said to be
	\begin{enumerate}
		\item (Lyapunov) stable, if for all~${\epsilon>0}$ there exists a scalar~${\delta_1\in(0,\epsilon)}$ such that~${\LInfNorm{y_0}+\norm{\mathbf{x}(0)}_{e_{1}}<\delta_1}$ implies~${\LInfNorm{y_k}<\epsilon}$, for all~${k\in\nnum}$,
		\item asymptotically stable, if it is Lyapunov stable and there exists~${\delta_2>0}$ such that~${\LInfNorm{y_0}+\norm{\mathbf{x}(0)}_{e_{1}}<\delta_2}$ and~${\mathbf{x}(0)\in c_0}$ implies~${\LInfNorm{y_k}\to 0}$,
		\item exponentially stable, if it is asymptotically stable, and there exist~${\delta_3>0}$ and continuous increasing functions~${K:(0,1)\to [1,\infty)},{\gamma:(0,1)\to(0,1)}$, such that~${\LInfNorm{y_0}+\norm{\mathbf{x}(0)}_{e_{\lambda}}<\delta_3}$ implies
		\begin{equation}
		\LInfNorm{y_k}\leq K(\lambda)\gamma(\lambda)^k(\LInfNorm{y_0}+\norm{\mathbf{x}(0)}_{e_{\lambda}}),
		\label{eq:expstab}
		\end{equation}
		 for all~${k\in\nnum}$ and~${\lambda\in(0,1)}$.
	\end{enumerate}
\end{defn}

In the rest of the paper, since the origin is the only equilibrium of interest, we will simply say that the DRP~\eqref{eq:nonlinear} is (Lyapunov)/asymptotically/exponentially stable. In addition, we will say that the DRP~\eqref{eq:nonlinear} is globally asymptotically~(exponentially) stable if~$\delta_2$~($\delta_2$ and~$\delta_3$) can be chosen to be arbitrarily large. A salient feature of the exponential stability definition above is the dependency of the performance on the convergence speed~$\lambda$ of~$\mathbf{x}(0)$, expressed via the functions~$K$ and~$\gamma$, which are continuous and increasing to be physically meaningful. In addition, since 0 is an equilibrium solution for~\eqref{eq:nonlinear1D}, which is Lipschitz with respect to~${(\chi(0),u)}$ by Lemma~\ref{lem:lipschitz}, it is straightforward to show that the stability notions above translate directly to the state trajectory.

We will also be considering the case~${\mathbf{x}(0)=0}$. We will refer to any such DRP as a \textit{zero initial states}~{(0-i.s.)} system or process. The~0-i.s. system will be defined to be Lyapunov, asymptotically, or exponentially stable if the notions defined above hold for the case of~${\mathbf{x}(0)=0}$; obviously the~0-i.s. system is (asymptotically/exponentially) stable if the actual system is (asymptotically/exponentially) stable. Note that~\eqref{eq:expstab} is necessary~\textit{and} sufficient for~0-i.s exponential stability.

\section{Stability of LTV Differential Processes}
\label{sec:linear}
In this section, we will focus on systems where~$f$ and~$g$ are linear with respect to their first two arguments for fixed~$t{\in[0,T]}$, and relax the continuous differentiability assumption to that of continuity;~i.e. we will look at LTV differential processes of the form
\begin{equation}
	\left\{
	\begin{aligned}
		\dot{x}_{k+1}(t)		&= A(t)x_{k+1}(t)+B(t)y_{k}(t),\\
				y_{k+1}(t)		&= C(t)x_{k+1}(t)+D(t)y_{k}(t),
	\end{aligned}
	\right.
\label{eq:LTV}
\end{equation}
for all~${(t,k)\in[0,T]\times\mathbb{N}}$, where~${A,B,C,D}$ are continuous real matrices of appropriate size.


\subsection{0-i.s Stability and the Spectral Radius}

Similar to the nonlinear case, given the LTV system described by the quadruple~${(A,B,C,D)}$, we denote by~$G_x$ the state response to the input and the initial condition, and by~$G_y$ the mapping from the input and the state to the output. The LTV operator~$G$ is defined so that
\[
	(w,\chi)=G(u,\chi(0))=(G_y(u,G_x(u,\chi(0))),G_x(u,\chi(0))),
\]
and the~0-i.s. output response~${G_0(.)\triangleq\pi_y(G(.,0))}$, where~$\pi_y$ is the standard projection onto~$\mathcal{Y}$. We will first consider the~0-i.s. system described by the discrete system~${y_{k+1}=G_0y_k}$ on~$\mathcal{Y}$. We have the following claim about~$G_0$:
\begin{claim}
	\label{claim:Lp}
	The operator~$G_0$ is bounded in~$\mathcal{L}_p$, for any~${p\in[1,\infty]}$.
\end{claim}

Claim~\ref{claim:Lp} makes intuitive sense since linear systems do not have finite escape time. The formal proof of this argument relies on the continuity of state matrices~(and hence that of the state-transition matrix) and~\cite[Theorem~75]{vidyasagar}; see Appendix~\ref{app:claimproof}. As such, we will expand the space~$\mathcal{Y}$ to~$\mathcal{L}_\infty$, and more generally~$\mathcal{L}_p$. The stability problem is relatively simple for linear systems as expected:~Exponential stability can be conveniently evaluated by the following spectral radius condition, which can easily be proven by Gelfand's spectral radius formula~${\rho(G_0)=\lim_{k\to\infty}\norm{G_0^k}_{\mathcal{L}_p}^{1/k}}$~\cite{weiss}:
\begin{thm}
The~0-i.s. linear system~\eqref{eq:LTV} is exponentially stable~(in~$\mathcal{L}_p$) if and only if~$\rho(G_0)<1$.
\end{thm}
\begin{rem}
	In general, the condition~${\rho(G_0)<1}$ is sufficient for asymptotic stability, whereas~${\rho(G_0)\leq 1}$ is necessary~\cite{przyluski}. This issue is circumvented in~[page~44]\cite{rogers} by requiring asymptotic stability to be a local property around a nominal operator.
\end{rem}


\subsection{Computation of the Spectral Radius}

The computation of the spectral radius will be similar to the procedure outlined for the time-invariant case in~\cite{rogers}. Let~${P_z(t)\triangleq zI-D(t)}$, where~${z\in\mathbb{C}}$. {\color{black}It is easy to see that the operator~$zI-G_0$ mapping~$u$ to~$\eta$, given by
\begin{equation*}
	\label{eq:singular}
	\left\{
	\begin{aligned}
		\dot{\chi}(t)		&= A(t)\chi(t)+B(t)u(t),\\
		\eta(t)					&=-C(t)\chi(t)+P_z(t)u(t),
	\end{aligned}
	\right.
\end{equation*}
for all~${t\in[0,T]}$, is invertible if~${|z|>\alpha\triangleq\sup_{t\in [0,T]}\rho(D(t))}$. In addition,~$(zI-G_0)^{-1}$ is bounded~(in~$\mathcal{L}_p$) by the bounded inverse theorem. Hence,~$\rho(G_0)\leq \alpha$.}

Otherwise, given any~${\epsilon>0}$, let~${z\in\mathbb{C}}$ be a number such that~$P_z(t)$ is singular for some~${t\in [0,T)}$ and~${|z|>\alpha-\epsilon}$. Such a~$z$ exists since the spectral radius of~$D$ varies continuously. Define~${s \triangleq \min\{t\in[0,T]: \mathrm{det}(P_z(t))=0\}}$, and set~${\eta(t)=\phi\mathbf{1}(t-s)}$, where~$\phi$ is orthogonal to the range of~$P_z(s)$, and~$\mathbf{1}(.)$ is the Heaviside step function. Assume that there exists a~$u\in\mathcal{L}_{\infty}$ that achieves~$\eta$ almost everywhere. Obviously, the input~$u=0$ and state~$\chi=0$, almost everywhere on~$[0,s)$. Define
\[
	\mu(t)\triangleq\norm{\phi-P_z(t)u(t)+C(t)\chi(t)}_2, \quad \forall t\in[s,T].
\]
By~\eqref{eq:singular},~$\mu=0$ almost everywhere on~$[s,T]$. Moreover, since~$\chi$ is continuous\footnote{See~\cite[page 48]{rugh} for piecewise continuous~$u$ and~\cite[Theorem~II.4.6]{warga} for integrable~$u$.} by~\eqref{eq:singular},~$\chi(s)=0$. Now let~$\Psi$ be an orthogonal projection matrix, onto the span of~$\phi$. {\color{black} Using the reverse triangle inequality, by orthogonality, it is easy to show
	\[
		\mu(t)\geq \norm{\phi}_2-\left(\norm{\Psi P_z(t)u(t)}_2+\norm{C(t)\chi(t)}_2\right),
	\]
	for all~$t\in[0,T]$.} Clearly,~$\mu(s)\geq\norm{\phi}_2$. In addition, since~$P_z,C,\chi$ are continuous,~$\chi(s)=0$, and~$\Psi P_z(s)=0$, the scalar~$\sup_{\tau\in[s,t)}\left(\norm{\Psi P_z(\tau)}_2+\norm{C(\tau)\chi(\tau)}_2\right)$ can be made arbitrarily small as~$t$ approaches~$s$ from the right. Consequently, given any~$u\in\mathcal{L}_{\infty}$, the essential supremum of~$\norm{C(\tau)\chi(\tau)}_2+\norm{\Psi P_z(\tau)u(\tau)}_2$ can be made arbitrarily small almost everywhere on~$[s,t)$ as~$t$ approaches~$s$ from the right. But then,~$\mu(t)\geq \varsigma>0$, almost everywhere on~$[s,t)$ for some~$t>s$ and constant~$\varsigma$, contradicting the fact that~$\mu=0$. It follows that~${zI-G_0}$ is not surjective. Therefore,~${\rho(G_0)=\alpha}$.


\subsection{Stability under Nonzero Initial States}\label{subsec:nonzero}

Let~${H(.)\triangleq \pi_y(G(0,.))}$ be the natural response of the LTV system to initial conditions. Then the solution of~\eqref{eq:LTV} can be given as
\begin{equation*}
	y_{k}=G_0^ky_0+\sum _{i=1}^k G_0^{k-i}Hx_{i}(0), \quad \forall k\in\nnum.
	\label{eq:solution}
\end{equation*}
Now if~${\rho(G_0)<1}$, by Gelfand's spectral radius formula, there exist scalars~${M>0}$ and~${\zeta\in (0,1)}$ such that~${\LInfNorm{G_0^k}\leq M\zeta^k}$ for all~${k\in\nnum}$. Therefore,
\begin{equation}
	\LInfNorm{y_{k}}\leq M\left(\zeta^k\LInfNorm{y_0}	+\LInfNorm{H}\sum _{i=1}^k \zeta^{k-i}\norm{x_{i}(0)}\right),
\label{eq:solutionbound}
\end{equation}
for all~${k\in\nnum}$, where~$H$ is bounded due to the finite-time assumption.\footnote{See the discussion of Claim~\ref{claim:Lp}.} {\color{black}When~${\norm{\mathbf{x}(0)}_{e_1}<\infty}$, it is easy to bound the right-hand side of~\eqref{eq:solutionbound} as a linear function of~$\left(\LInfNorm{y_0}+\norm{\mathbf{x}(0)}_{e_1}\right)$.} Therefore, the LTV system is stable. Now assume in addition that~${\mathbf{x}(0)\in c_0}$, and consider the partial sum in the second term of the right hand side of~\eqref{eq:solutionbound},~${S_k\triangleq\sum _{i=1}^k \zeta^{k-i}\norm{x_{i}(0)}\geq 0}$,~for all~${k\in\nnum}$. Then, it is easy to verify~${S_{k+1}=\zeta S_k+\norm{x_{k+1}(0)}\geq 0}$ for all~${k\in\nnum}$, so by Claim~\ref{claim:asymp},~${S_k\to 0}$. Therefore, we can conclude by~\eqref{eq:solutionbound} that~${y_{k}\to 0}$ if~${\mathbf{x}(0)\in c_0}$ and~${\rho(G)<1}$.

Finally, consider the case~${\mathbf{x}(0)\in e_{\lambda}}$ for some~${\lambda\in(0,1)}$. From~\eqref{eq:solutionbound}
\begin{equation*}
	\begin{aligned}
	\norm{y_{k}}&\leq M\left(\zeta^k\LInfNorm{y_0}+\LInfNorm{H}\norm{\mathbf{x}(0)}_{e_{\lambda}}\sum _{i=1}^k \zeta^{k-i}\lambda^{i-1}\right)\\
		&\leq M\left(\zeta^k\LInfNorm{y_0}+\LInfNorm{H}\norm{\mathbf{x}(0)}_{e_{\lambda}}k\bar{\lambda}^{k-1}\right),
	\end{aligned}
\end{equation*}
where~${\bar{\lambda}\triangleq\max\{\zeta,\lambda\}}$, so by Claim~\ref{claim:exp}
\begin{multline*}
	\norm{y_{k}}\leq M\zeta^k\LInfNorm{y_0}\\
	+M\LInfNorm{H}\norm{\mathbf{x}(0)}_{e_{\lambda}}\frac{2}{1-\bar{\lambda}}\left(\frac{1+\bar{\lambda}}{2}\right)^{k},
\end{multline*}
and since~${\zeta\leq \bar{\lambda}< (1+\bar{\lambda})/2<1}$,
\begin{multline}\label{eq:Kg}
\LInfNorm{y_{k}} \leq \overbrace{M\max\left\{1,\frac{2\LInfNorm{H}}{1-\bar{\lambda}}\right\}}^{K_G(\bar{\lambda})} \\
\times\left(\underbrace{(1+\bar{\lambda})/2}_{\gamma_G(\bar{\lambda})}\right)^{k}(\LInfNorm{y_0}+\norm{\mathbf{x}(0)}_{e_{\lambda}}), \quad \forall k\in\nnum.
\end{multline}
Noting that~${K_G(\max\{\zeta,\lambda\})}$ and~${\gamma_G(\max\{\zeta,\lambda\})}$ defined in~\eqref{eq:Kg} are both continuous and increasing in~$\lambda$ on~${(0,1)}$, we can conclude the system to be exponentially stable. With this, our findings can be summarized as follows:
\begin{thm}\label{thm:G}
	For the LTV DRP~\eqref{eq:LTV}, the following are equivalent:
	\begin{enumerate}
		\item The DRP~\eqref{eq:LTV} is globally exponentially stable.
		\item The 0-i.s. DRP~\eqref{eq:LTV} is globally exponentially stable.
		\item The condition~$max_{t\in [0,T]}\rho(D(t))<1$ holds.
	\end{enumerate}
\end{thm}
\begin{rem}\label{rem:Lp}
The analysis of Section~\ref{subsec:nonzero} extends to any~$\mathcal{L}_p$ norm since~$\rho(G_0)\leq\alpha$ for all~$p\in[1,\infty]$. Therefore,~${\alpha<1}$ implies global exponential stability in~$\mathcal{L}_p$.
\end{rem}

\section{Linearized Stability of DRPs}
\label{sec:linearized}
We will now establish the equivalence between exponential stability of a nonlinear DRP of the form~\eqref{eq:nonlinear} with that of its linearization. The linearization of~\eqref{eq:nonlinear} will mirror that of the~1D case, in other words, we will be linearizing the differential operator~\eqref{eq:nonlinear1D} as is typical in feedback control. This will be done as follows: Since~$f$ and~$g$ are continuously differentiable,
\begin{equation}
	\left\{
	\begin{aligned}
	\dot{\chi}(t) &=\bar{A}(t)\chi(t) +\bar{B}(t)u(t) + b(\chi(t),u(t),t),\\
	w(t)					&=\bar{C}(t)\chi(t) +\bar{D}(t)u(t) + d(\chi(t),u(t),t),
	\end{aligned}
	\right.
\label{eq:linearx}
\end{equation}
for some continuous functions~$b$ and~$d$, as
\begin{equation*}
	\begin{aligned}
		&\bar{A}(t) \triangleq\frac{\partial f}{\partial \chi}(0,0,t), \quad &\bar{B}(t) \triangleq\frac{\partial f}{\partial u}(0,0,t),\\
		&\bar{C}(t) \triangleq\frac{\partial g}{\partial \chi}(0,0,t), \quad &\bar{D}(t) \triangleq\frac{\partial g}{\partial u}(0,0,t),
	\end{aligned}
\end{equation*}
are continuous. Consequently, the linearization of~\eqref{eq:nonlinear} will be defined as the following~2D system:
\begin{equation}
	\left\{
	\begin{aligned}
		\dot{\bar{x}}_{k+1}(t)&=\bar{A}(t)\bar{x}_{k+1}(t) +\bar{B}(t)\bar{y}_k(t),\\
		\bar{y}_{k+1}(t)			&=\bar{C}(t)\bar{x}_{k+1}(t) +\bar{D}(t)\bar{y}_k(t),
	\end{aligned}
	\right.
	\label{eq:linearized}
\end{equation}
for all~${(t,k)\in[0,T]\times\nnum}$, with boundary conditions satisfying~${\mathbf{\bar{x}}(0)=\mathbf{x}(0)}$ and~${\bar{y}_0=y_0}$.


\subsection{Asymptotics of the Nonlinear Perturbations}

Let~$f_i$ be the~$i$-th output of~$f$. Since~$f$ is continuously differentiable in~$Z$ and~${f(0,0,t)=0}$, by the multivariable mean value theorem, there exists a point ${(\xi_i^*,\upsilon_i^*)}$ on the line segment connecting~${(\xi,\upsilon)}$ to the origin such that
\[
	f_i(\xi,\upsilon,t)=
	\begin{bmatrix}
		\frac{\partial f_i}{\partial \xi}(\xi_i^*,\upsilon_i^*,t) & \frac{\partial f_i}{\partial \upsilon}(\xi_i^*,\upsilon_i^*,t)
	\end{bmatrix}
	\begin{bmatrix}
		\xi \\ \upsilon
	\end{bmatrix}
\]
in a neighborhood of~${0\in\real^n\times\real^m}$. Equivalently,
\begin{multline*}
	f_i(\xi,\upsilon,t)=
	\begin{bmatrix}
		\bar{A}_i(t) & \bar{B}_i(t)
	\end{bmatrix}
	\begin{bmatrix}
		\xi \\ \upsilon
	\end{bmatrix}\\
	+
	\underbrace{
	\begin{bmatrix}
		\left(\frac{\partial f_i}{\partial \xi}(\xi_i^*,\upsilon_i^*,t)-\bar{A}_i(t)\right) & \left(\frac{\partial f_i}{\partial \upsilon}(\xi_i^*,\upsilon_i^*,t)-\bar{B}_i(t)\right)
	\end{bmatrix}
	\begin{bmatrix}
		\xi \\ \upsilon
	\end{bmatrix}}_{b_i(\xi,\upsilon,t)},
\end{multline*}
where~$\bar{A}_i$ and~$\bar{B}_i$ are the~$i$-th rows of~$\bar{A}$ and~$\bar{B}$, respectively, and~$b_i$ is the~$i$-th output of~$b$. Now let~${Q_i\triangleq \partial f_i/\partial \xi}$. The function~$Q_i$ is continuous in~$Z$ because~$f$ is continuously differentiable in~$Z$. Hence, by the Heine-Cantor theorem,~$Q_i$ is uniformly continuous in~$Z$. Therefore, for all~${\epsilon>0}$ there exists~${\delta_o>0}$ such that
\[
	\norm{(\xi,\upsilon)}<\delta_o \implies \norm{Q_i(\xi,\upsilon,t)-\bar{A}_i(0,0,t)}<\epsilon,
\]
for every~${(\xi,\upsilon,t)\in Z}$, since~${Q_i(0,0,t)=\bar{A}_i(0,0,t)}$. Using similar arguments for~${\partial f_i/\partial \upsilon,\partial g_i/\partial \xi,\partial g_i/\partial \upsilon}$, we can conclude that for all~${\epsilon>0}$ there exists~${\delta_O>0}$ satisfying
\begin{multline}
\label{eq:O(x)}
	\norm{(\xi,\upsilon)}<\delta_O\\
	\implies \norm{(b(\xi,\upsilon,t),d(\xi,\upsilon,t)}<\epsilon \norm{(\xi,\upsilon)},
\end{multline}
for every~${(\xi,\upsilon,t)\in Z}$.

\subsection{$\mathcal{L}_{\infty}$ Asymptotics of the Linearization Error}

Next, let us consider the LTV system defined by the matrices~${\bar{A},\bar{B},\bar{C},\bar{D}}$:
\begin{equation}
	\left\{
	\begin{aligned}
	\dot{\bar{\chi}}(t)	&=\bar{A}(t)\bar{\chi}(t) +\bar{B}(t)\bar{u}(t),\\
	\bar{w}(t)					&=\bar{C}(t)\bar{\chi}(t) +\bar{D}(t)\bar{u}(t),
	\end{aligned}
	\right.
	\label{eq:linearized1D}
\end{equation}
for all~${t\in[0,T]}$, where~${\bar{\chi}(0)=\chi(0)}$. The~0-i.s. input-output response~$\bar{G}_0$ and the initial state response~$\bar{H}$ will be defined for this system as in Section~\ref{sec:linear}. Subtracting~\eqref{eq:linearized1D} from~\eqref{eq:linearx},
\begin{equation*}
	\left\{
	\begin{aligned}
	\dot{\tilde{\chi}}(t)	&=\bar{A}(t)\tilde{\chi}(t) + \bar{B}(t)\tilde{u}(t)+ b(\chi(t),u(t),t),\\
	\tilde{w}(t)					&=\bar{C}(t)\tilde{\chi}(t) + \bar{D}(t)\tilde{u}(t)+d(\chi(t),u(t),t),
	\end{aligned}
	\right.
\end{equation*}
where~${\tilde{\chi}(t)\triangleq{\chi}(t)-\bar{\chi}(t)}$,~${\tilde{w}(t)\triangleq{w}(t)-\bar{w}(t)}$, and similarly~${\tilde{u}(t)\triangleq{u}(t)-\bar{u}(t)}$. Define the mapping~$\varphi$ so that
\[
	(\varphi(\chi,u))(t)=(b(\chi(t),u(t),t),d(\chi(t),u(t),t)).
\]
Then the output error~$\tilde{w}$ is given by
\begin{equation}
\tilde{w}=\bar{G}_0\tilde{u}+\Omega(\varphi(\chi,u)),
\label{eq:operator}
\end{equation}
where~$\Omega$ represents the~$\mathcal{L}_{\infty}$ stable input-output response of an LTV system with state matrices~${A,[\begin{matrix}I & 0\end{matrix}],C,[\begin{matrix}0 & I \end{matrix}]}$. The following lemma will define the asymptotic behavior of~$\varphi$ with respect to~${(u,\chi(0))}$; see Appendix~\ref{app:proofs} for a proof.
\begin{lem}
\label{lem:O(x)}
	For all~$\epsilon>0$, there exists~$\delta^*\in(0,\epsilon)$ such that
	\begin{multline*}
		\LInfNorm{u}+\norm{\chi(0)}<\delta^*\\
		\implies \LInfNorm{\varphi(\chi,u)}\leq\epsilon(\LInfNorm{u}+\norm{\chi(0)}).
	\end{multline*}
\end{lem}


\subsection{Necessary and Sufficient Conditions for Exponential Stability}

We first assume that the~0-i.s. linearized system is exponentially stable so that~${\LInfNorm{\bar{G}_0^k}\leq \bar{M}\bar{\zeta}^k}$ for all~${k\in\nnum}$, for some~${\bar{M}\geq1,\bar{\zeta}\in(0,1)}$. With this, let~${N\in\nnum}$ such that~${\bar{M}\bar{\zeta}^N<1}$. We will need the subsequent result, which follows easily from Lipschitz continuity of~${\Gamma}$~(Lemma~\ref{lem:lipschitz}):
\begin{lem}
\label{lem:finiteL}
	There exist scalars~$\delta_{\rm fh}>0$ and~$L_{\rm fh}\geq 1$ so that~$\LInfNorm{y_0}+\norm{\mathbf{x}(0)}_{e_1}<\delta_{\rm fh}$ implies
	\[
		\LInfNorm{y_k}<L_{\rm fh}(\LInfNorm{y_0}+\norm{\mathbf{x}(0)}_{e_1}),
	\]
for all~${k\in\{0,1,\dots,N-1\}}$.
\end{lem}

\begin{prop}\label{prop:sufficiency}
	The nonlinear system~\eqref{eq:nonlinear} is exponentially stable if its linearization~\eqref{eq:linearized} is exponentially stable.
\end{prop}

\begin{figure}[tbhp]
	\centering
	\includegraphics[width=1\columnwidth]{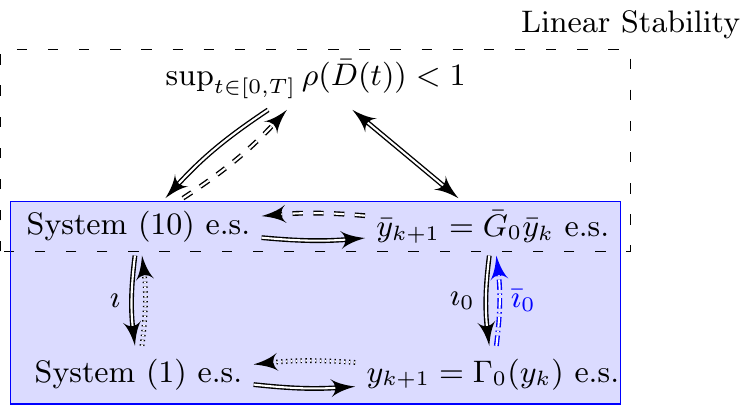}
	\caption{Implication diagram for exponential stability~(e.s.): The linear exponential stability diagram was stated in Theorem~\ref{thm:G}, where the dashed implication arrows were established by proving the solid ones. For the nonlinear case, implications~$\imath,\imath_0$ are proven in Proposition~\ref{prop:sufficiency}. Proving implication~$\bar{\imath}_0$ will close the loop and allow us to conclude the dotted implication arrows.}
	\label{fig:implication}
\end{figure}

The proof of this proposition is rather involved and as such given in Appendix~\ref{app:sufficiency} for a more compact presentation. To establish the converse of this result, we will follow an indirect route that is much easier compared to a direct proof. Specifically, we will show that nonlinear exponential stability implies linear exponential stability for the~0-i.s case. Since the~0-i.s. systems given by the operators~${\Gamma_0(.)\triangleq\pi_y(\Gamma(.,0))}$ and~$G_0$ are in essence discrete systems evolving on~$\mathcal{Y}$, we will be relying on the following forward Lyapunov-like theorem. This will allow us to finalize our main result by aid of Theorem~\ref{thm:G}, as can be seen in Fig.~\ref{fig:implication}. For obvious reasons, a functional satisfying the conditions of Theorem~\ref{thm:expstable} will be called a Lyapunov functional.
\begin{thm}\label{thm:expstable}
	Let~$F_0\in\{\Gamma_0,G_0\}$. Then, the 0-i.s. DRP given by the recursion~${y_{k+1}=F_0(y_k)}$ on~$\mathcal{Y}$ is exponentially stable if and only if there exist a functional~${V:\mathcal{Y}\to\real}$ and positive scalars~${c_1,c_2,c_3}$, with~${c_2>c_3}$, such that~${c_1\LInfNorm{y}\leq V(y)\leq c_2\LInfNorm{y}}$ and~${V(F_0(y))-V(y)\leq -c_3\LInfNorm{y}}$~ in a neighborhood of the origin.
\end{thm}
\begin{pf}
Sufficiency is obvious and is therefore omitted. The necessity part can be proven by construction as follows: Assume that the system is exponentially stable, then there exist~${K>1}$,~${\delta_3>0}$ and~${\gamma\in[0,1)}$ so that~${\LInfNorm{F_0^k(y)}\leq K\gamma^k\LInfNorm{y}}$ holds for all~${y\in\mathcal{Y}}$ with~${\norm{y}<\delta_3}$. Let~$N$ be an integer so~${K\gamma^N<1}$. Then, it is easy to show~${V(y)\triangleq\sum_{i=0}^{N-1}\LInfNorm{F_0^i(y)}\geq\LInfNorm{y}}$ satisfies the conditions of the theorem for all~${y\in\mathcal{Y}}$ with~${\LInfNorm{y}<\delta_3}$.
\end{pf}
\begin{prop}
\label{prop:necessity}
	The linearization~\eqref{eq:linearized} of the nonlinear system is~0-i.s. exponentially stable if the nonlinear system~\eqref{eq:nonlinear} is~0-i.s.exponentially stable.
\end{prop}
\begin{pf}
	Let~$F_0=\Gamma_0$, and let~$V$ be the Lyapunov functional from the proof of Theorem~\ref{thm:expstable}. Then,~$V$ satisfies~${c_1\LInfNorm{y}\leq V(y)\leq c_2\LInfNorm{y}}$, and the difference of~$V$ with respect to the linear operator~$\bar{G}_0$ is
	\begin{equation*}
		\begin{aligned}
			\Delta V(y)	&\triangleq V(\bar{G}_0y)-V(y)\\
									&=(V(\bar{G}_0y)-V(\Gamma_0(y)))+(V(\Gamma_0(y))-V(y))\\
									&\leq(V(\bar{G}_0y)-V(\Gamma_0(y)))-c_3\LInfNorm{y},
		\end{aligned}
	\end{equation*}
	around the origin for some positive~${c_1,c_2,c_3}$ with~${c_2>c_3}$, as the nonlinear DRP is exponentially stable. The functional~$V$ is locally Lipschitz because it is a sum of locally Lipschitz functionals;~${\Gamma_0^i}$ is locally Lipschitz for any~${i\in\nnum}$ by Lemma~\ref{lem:lipschitz}. Furthermore, from~\eqref{eq:operator},
	\[
		\Bar{G}_0y-\Gamma_0(y)=\Omega(\varphi(\Gamma_x(y,0),y)).
	\]
	Recall that~$\Omega$ is~$\mathcal{L}_{\infty}$ stable, and~$\Gamma_x$ is locally Lipschitz. Hence, for any~${\epsilon>0}$, by Lemma~\ref{lem:O(x)}, there exists~${\delta^*>0}$ so~${\LInfNorm{y}<\delta^*}$ implies~${|V(\Gamma_0(y))-V(\bar{G}_0y)|\leq\epsilon\LInfNorm{y}}$, and therefore for any~${\bar{c}_3\in(0,c_3)}$, there exists a~${\bar{\delta}_3>0}$ so that~${\LInfNorm{y}\leq\bar{\delta}_3}$ implies
	\[
	\Delta V(y)\leq(V(\bar{G}_0y)-V(\Gamma_0(y)))-c_3\LInfNorm{y}\leq \bar{c}_3\LInfNorm{y}.
	\]
	By Theorem~\ref{thm:expstable}, it follows that the linearization~\eqref{eq:linearized} is~0-i.s. exponentially stable.
\end{pf}

We are now ready to state our main result, which summarizes the findings of Theorem~\ref{thm:G} and Propositions~\ref{prop:sufficiency} and~\ref{prop:necessity} as given below:
\begin{thm}\label{thm:main}
	For the nonlinear DRP~\eqref{eq:nonlinear} and its linearization~\eqref{eq:linearized}, the following are equivalent:
	\begin{enumerate}
		\item The DRP~\eqref{eq:nonlinear} is exponentially stable.
		\item The~0-.i.s. DRP~\eqref{eq:nonlinear} is exponentially stable.
		\item The DRP~\eqref{eq:linearized} is globally exponentially stable.
		\item The~0-.i.s. DRP~\eqref{eq:linearized} is globally exponentially stable.
		\item The condition~$\max_{t\in[0,T]}\rho(\bar{D}(t))<1$ holds.
	\end{enumerate}
\end{thm}

\section{Applications: Picard Iterations and ILC}
\label{sec:picardILC}

{
We now present two applications of Theorem~\ref{thm:main}.


\subsection{Picard Iterates with Varying Initial Conditions}

The Picard-Lindel\"{o}f theorem guarantees the existence and uniqueness of the solution~$x^*$ of the differential equation~${\dot{x}(t)=f(x(t),t)}$ with initial condition~${x(0)=x^*_0}$ for small~${T}$. The existence of this solution is proven by a recursive process, whose convergence is shown by the contraction mapping theorem. These iterates can be expressed as the DRP
\begin{equation*}
	\left\{
	\begin{aligned}
		\dot{x}_{k+1}(t)		&= f(y_{k}(t),t), \quad x_{k+1}(0)=x^*_0,\\
				y_{k+1}(t)		&= x_{k+1}(t),
	\end{aligned}
	\right.
\end{equation*}
for all~${(t,k)\in[0,T]\times\nnum}$. The time-varying transformation~${(x_k(t),y_k(t))\mapsto(x_k(t)-x^*(t),y_k(t)-x^*(t))}$ translates the equilibrium to~0, uniformly in time:
\begin{equation*}
	\left\{
	\begin{aligned}
		\dot{\underline{x}}_{k+1}(t)		&= \underline{f}(\underline{y}_{k}(t),t), \quad \underline{x}_{k+1}(0)=0,\\
				\underline{y}_{k+1}(t)		&= \underline{x}_{k+1}(t),
	\end{aligned}
	\right.
\end{equation*}
with~${\underline{f}(\chi,t)\triangleq f(\chi+x^*(t),t)-\dot{x}^*(t)}$, for all~${t\in[0,T]}$ and~${k\in\nnum}$. This resulting system satisfies continuous differentiability assumptions around the new equilibrium since the fixed point~$x^*$ is \textit{twice} continuously differentiable by virtue of~$f$ being continuously differentiable. Now, we can conclude that Picard iterates form an exponentially stable DRP when~${y_0-x^*}$ and~${\mathbf{x}(0)-x^*_0}$ are small enough. Hence, the iterates converge to~$x^*$ for every~$\mathbf{x}(0)$ with ${\mathbf{x}(0)-x^*_0\in c_0}$, e.g. for \textit{nonconstant initial state sequences that converge to~$x_0$}, when the boundaries are close to the equilibrium.


\subsection{ILC with Static Nonlinear Update Laws}

The second application of Theorem~\ref{thm:main} addresses the ILC problem of iteratively constructing the feedforward input~$u^*$ given a desired output~$y_{\mathrm{des}}$ so that
}
\begin{equation*}
	\left\{
	\begin{aligned}
		\dot{x}^*(t)		&= f(x^*(t),u^*(t),t),\\
				y_{\mathrm{des}}(t)		&= g(x^*(t),u^*(t),t),
	\end{aligned}
	\right.
\end{equation*}
for all~${t\in[0,T]}$. We consider the ILC system, where the continuously differentiable function~$l$ satisfies~${l(0,t)=0}$,
\begin{equation*}
	\left\{
	\begin{aligned}
		\dot{x}_{k+1}(t)		&= f(x_{k+1}(t),u_{k+1}(t),t),\\
				y_{k+1}(t)		&= g(x_{k+1}(t),u_{k+1}(t),t),\\
				u_{k+1}(t)		&= u_{k}(t)+l(e_k(t),t),
	\end{aligned}
	\right.
\end{equation*}
and~${e_k\triangleq y_k-y_{\mathrm{des}}}$, for all~${(t,k)\in[0,T]\times\nnum}$. This static~(in-time) update law is based on the internal model principle in the iteration domain, and guarantees perfect tracking in the limit for all achievable~$y_{\mathrm{des}}$ when stable. Following a transformation akin to the one for Picard iterates, we can rewrite the system as
\begin{equation*}
	\left\{
	\begin{aligned}
		\dot{\underline{x}}_{k+1}(t)		&= \underline{f}(\underline{x}_{k+1}(t),\underline{u}_{k}(t),e_k(t),t),\\
				\begin{bmatrix}
				e_{k+1}(t)\\ \underline{u}_{k+1}(t)
				\end{bmatrix} &=
				\begin{bmatrix}
				\underline{g}(\underline{x}_{k+1}(t),\underline{u}_{k}(t),e_k(t),t),\\
				\underline{u}_{k}(t)+l(e_k(t),t)
				\end{bmatrix},
	\end{aligned}
	\right.
\end{equation*}
with
\[
	\underline{g}(\chi,u,\theta,t)\triangleq g(\chi+x^*(t),u+u^*(t)+l(\theta,t),t)-y_{\mathrm{des}}(t),
\]
for all~${(t,k)\in[0,T]\times\nnum}$. Observe that~$e_0$ depends on~$\underline{u}_0$, so~${(e_0,\underline{u}_0)}$ cannot be arbitrarily chosen, and thus it is difficult to derive necessary stability conditions. Nevertheless, letting
\[
	\underline{D}(t)\triangleq \frac{\partial g}{\partial u}(x^*(t),u^*(t),t), \quad \underline{L}(t)\triangleq \frac{\partial l}{\partial \theta}(0,t),
\]
 for all~${t\in[0,T]}$, the system is exponentially stable if
	\[
	\max_{t\in[0,T]}\rho\left(
	\begin{bmatrix} \underline{D}(t) \\ I \end{bmatrix}\begin{bmatrix} \underline{L}(t) & I\end{bmatrix}\right) =\max_{t\in[0,T]}\rho(I+\underline{L}(t)\underline{D}(t))<1,
\]
where the equality can be verified via simple eigenvector manipulations, with the equivalent condition being~${\max_{t\in[0,T]}\rho(I+\underline{D}(t)\underline{L}(t))<1}$ for square systems. Note that the same methodology can be used to derive spectral stability conditions with~$Q$ filtering; i.e. the update is of the form~${u_{k+1}(t)=\underline{Q}(t)u_{k}(t)+l(e_k(t),t)}$, a known robust stabilization factor in ILC algorithms.

The stability result derived above is the first eigenvalue based condition in the nonlinear ILC literature. Its significance further stems from the fact that it unifies several important results, such as continuous dependence of the tracking error on initial condition errors~\cite{heinzinger}, and the principle that the error term in the function~$l$ must be replaced with its~$\bar{n}$-th derivative for a relative degree~$\bar{n}$ system~\cite{ahn93}. Furthermore, it is among the first studies of ILC from a local perspective, which enables nonlinear time-varying update laws to be considered without resorting to saturation~\cite{tan}, and provides a rigorous basis to linearization in the context of ILC~\cite{bristow}.

\section{Illustrative Example}
\label{sec:example}
\begin{figure}[tbhp]
	\centering
	\includegraphics[width=1\columnwidth]{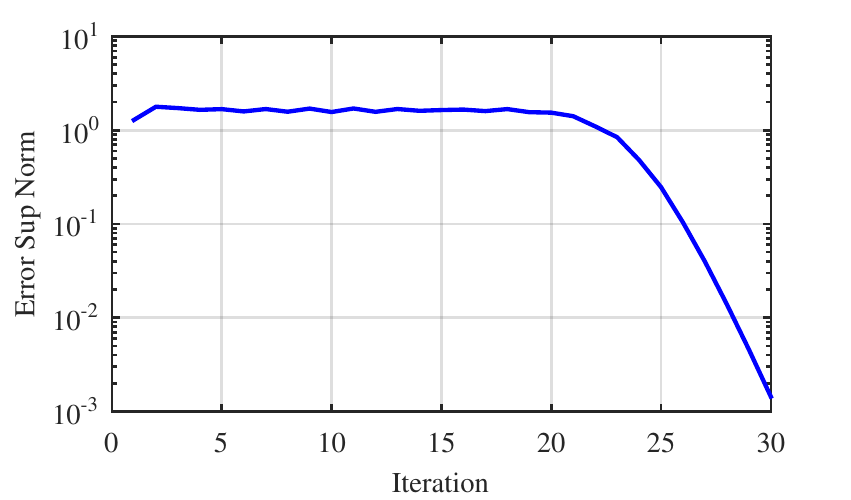}
	\caption{Evolution of~$\LInfNorm{e_k}$.}
	\label{fig:learning}
\end{figure}

Consider the actuated Van der Pol oscillator in normal form with a time-varying damping coefficient:
\begin{equation*}
	\left\{
	\begin{aligned}
		\dot{q}_1(t)&= q_2(t),\\
		\dot{q}_2(t)&= -q_1(t)+ \Xi(t)(1-(q_1(t))^2)q_2(t)+u(t),\\
		y						&= q_1(t),
	\end{aligned}
	\right.
\end{equation*}
where the damping coefficient~${\Xi(t)>0}$, and~${t\in[0,2]}$. The unforced oscillator is well-known to have an unstable equilibrium at the origin for all constant~${\Xi(t)>0}$. Our objective is to find an ILC update law in order to track the reference~${y_{\mathrm{des}}(t)=0.1\cos(2\pi t)}$. Since the relative degree is~2, we consider the update
\[
	u_{k+1}(t)=u_k(t)-(\ddot{y}_k(t)-\ddot{y}_{\mathrm{des}}(t)), \quad \forall (t,k)\in[0,2]\times\nnum.
\]
Then, it is easy to verify that this update law is stable since~${\ddot{y}(t)=\dot{q}_2(t)}$ and~${(\partial{\dot{q}_2}/\partial{u})(t)=1}$. Indeed, for~${\Xi(t)=4+0.5\sin(2\pi(10t))}$, Fig.~\ref{fig:learning} shows that the tracking error is exponentially decreased when~${u_0=0}$ and the initial conditions are randomly chosen to exponentially converge to~${(y_{\mathrm{des}}(0),\dot{y}_{\mathrm{des}}(0))=(0.1,0)}$ with convergence rate~$\lambda$~(also randomly chosen) and~$e_{\lambda}$ norm less than~$0.1$, \textit{without any stabilizing feedback}.


\section{Conclusion}
\label{sec:conclusion}

This paper addressed the problem of finding necessary and sufficient exponential stability conditions for a class of nonlinear repetitive processes and showed that a DRP is exponentially stable if and only if the state matrix~$\bar{D}$ of its linearization is uniformly Schur over the time interval~${[0,T]}$. To our knowledge, the work presented here is the first systematic study of local stability for nonlinear repetive processes. The findings of the paper are especially important since local stability is the precursor to global stability. The comprehensiveness of these results are reflected in the fact that they tie in the various existing results from nonlinear ILC analysis via a single framework. We hope that the analysis presented in the paper will pave the way for further research on nonlinear repetitive processes and other 2D systems, such as extensions to different classes of systems and the corresponding control strategies.


\appendix


\section{Proofs of Technical Results}
\label{app:proofs}
\begin{pf*}{Proof of Lemma~\ref{lem:lipschitz}}
	We begin by defining the set
\[
	\bar{\mathcal{Y}}\triangleq\{u\in\mathcal{Y}: u(t)\in Y, \forall t\in[0,T]\},
\]
and note that for any~${u\in\bar{\mathcal{Y}}}$,~${\bar{f}(\xi,t)\triangleq f(\xi,u(t),t)}$ is continuous in~${t\in[0,T]}$ for all~${\xi\in X}$ since~$f$ is continuous in~$Z$. Moreover, as~$f$ is continuously differentiable, it is also Lipschitz on the compact set~$Z$. That is, there exists a constant~$L_f$ such that
	\begin{multline*}
		\norm{f(\xi_1,\upsilon_1,\tau_1)-f(\xi_2,\upsilon_2,\tau_2)}\\
		\leq L_f\norm{(\xi_1,\upsilon_1,\tau_1)-(\xi_2,\upsilon_2,\tau_2)},
	\end{multline*}
	for all~${(\xi_1,\upsilon_1,\tau_1)}$ and~${(\xi_2,\upsilon_2,\tau_2)}$ in~$Z$. In turn, this implies that~${\norm{\bar{f}(\xi_1,t)-\bar{f}(\xi_2,t)}\leq L_f\norm{\xi_1-\xi_2}}$, for all~${\xi_1,\xi_2\in X}$ and~${t\in[0,T]}$, \textit{and} any~${u\in\bar{\mathcal{Y}}}$, so~${\bar{f}(\xi,t)}$ is Lipschitz with respect to~$\xi$, uniformly over time and the space of inputs. Now consider~${\dot{\chi_i}(t)=f({\chi_i}(t),u_i(t),t)}$, where the initial conditions and inputs satisfy the inequality~${\norm{\chi_i(0)}+\LInfNorm{u_i}<\delta}$, for each~${i\in\{1,2\}}$. By Assumption~\ref{assum:nonlinear}, the integral curves of both systems reside in~$X$. In addition,~${\bar{f}_i(\xi,t)\triangleq f(\xi,u_i(t),t)}$ is continuous in~${t\in[0,T]}$ for all~${\xi\in X}$, and Lipschitz with respect to~$\xi$ on~${X\times[0,T]}$, for each~${i\in\{1,2\}}$. Define the function~${\tilde{f}(\xi,t)\triangleq\bar{f}_2(\xi,t)-\bar{f}_1(\xi,t)}$, and rewrite the two systems as
	\begin{equation}
		\label{eq:khalil3_4}
		\begin{aligned}
			\dot{\chi}_1(t)&=\bar{f}_1({\chi_1}(t),t),\\
			\dot{\chi}_2(t)&=\bar{f}_1({\chi_2}(t),t)+\tilde{f}(\chi_2(t),t).
		\end{aligned}
	\end{equation}
	Since~$f$ is Lipschitz on~$Z$, as in the previous case where we showed that~$\bar{f}$ is Lipschitz with respect to its first argument, it follows that~${\norm{\tilde{f}(\xi,t)}\leq L_f\norm{u_1(t)-u_2(t)}}$ for all~${(\xi,t)\in X\times[0,T]}$. As~${u_1,u_2\in\bar{\mathcal{Y}}}$, this also means that~${\norm{\tilde{f}(\xi,t)}\leq L_f\LInfNorm{u_1-u_2}<R}$ for some~$R$, for all~${(\xi,t)\in X\times[0,T]}$ and~${u_1,u_2\in\bar{\mathcal{Y}}}$, since~$Y$ is compact. Now,~\eqref{eq:khalil3_4} satisfies all assumptions of Theorem~3.4 of~\cite{khalil}, which states that
	\begin{multline*}
		\norm{\chi_1(t)-\chi_2(t)}\leq \norm{\chi_1(0)-\chi_2(0)}e^{L_ft}\\
		+\LInfNorm{u_1-u_2}(e^{L_ft}-1), \quad \forall t\in[0,T],
	\end{multline*}
	therefore letting~${L_1=e^{L_fT}>1}$, we get
	\begin{multline*}
		\LInfNorm{\chi_1-\chi_2}\leq L_1(\norm{\chi_1(0)-\chi_2(0)}+\LInfNorm{u_1-u_2}),
	\end{multline*}
	when~${\norm{\chi_i(0)}+\LInfNorm{u_i}<\delta}$, for each~${i\in\{1,2\}}$.
	
	Continuous dependence of~$\omega$ on the pair~${(\chi(0),u)}$ can be shown in a similar way using continuous differentiability of~$g$; hence there exists~$L
	_2>0$ and~$\bar{\delta}\in(0,\delta)$ such that
	\begin{multline*}
		\LInfNorm{\omega_1-\omega_2}\leq L_2(\norm{\chi_1(0)-\chi_2(0)}+\LInfNorm{u_1-u_2})
	\end{multline*}
	when~${\norm{\chi_i(0)}+\LInfNorm{u_i}<\bar{\delta}}$, for each~${i\in\{1,2\}}$. Letting~$L=\max\{L_1,L_2\}>1$, the proof is complete.
\end{pf*}



\begin{pf*}{Proof of Lemma~\ref{lem:O(x)}}
	By Lemma~\ref{lem:lipschitz}, since the Lipschitz constant~${L\geq 1}$, the following is true:
	\begin{multline}
		\label{eq:dmy}
		\LInfNorm{u}+\norm{\chi(0)}<\bar{\delta}\\
		\implies \LInfNorm{(\chi,u)}\leq L(\LInfNorm{u}+\norm{\chi(0)}).
	\end{multline}
Moreover, by~\eqref{eq:O(x)}, for any~${\epsilon>0}$ there exists~${\delta_O^L>0}$ so
\[
	~{\LInfNorm{(\chi,u)}<\delta_O^L} \implies \LInfNorm{\varphi(\chi,u)}< (\epsilon/L)\LInfNorm{(\chi,u)}.
\]
Therefore, if~${\LInfNorm{u}+\norm{\chi(0)}<\delta^*<\min\{\bar{\delta},\delta_O^L/L,\epsilon\}}$, where~$\delta^*>0$ is arbitrary, it follows by~\eqref{eq:dmy} that
\[
{\LInfNorm{(\chi,u)}\leq L(\LInfNorm{u}+\norm{\chi(0)})<\delta_O},
\]
and consequently
\[
	\LInfNorm{\varphi(\chi,u)}< (\epsilon/L)\LInfNorm{(\chi,u)}<\epsilon(\LInfNorm{u}+\norm{\chi(0)}).
\]
\end{pf*}

\section{Discussion of Claim~\ref{claim:Lp}}
\label{app:claimproof}
We begin by noting that the matrices~${B,C,D}$ defining the operator~$G_0$ are continuous and hence bounded on~${[0,T]}$. Therefore, it is a straightforward matter to show that the multiplication operators defined by these matrices are bounded with respect to any~$\mathcal{L}_p$ norm,~${p\in[1,\infty]}$, so it will suffice to show that the time-varying convolution operator defined by the corresponding state-transition matrix is bounded. Because~$A$ is continuous, the state-transition matrix~$\Phi$ is continuously differentiable with respect to its first and second arguments on~${[0,T]^2}$~(see~\cite[page~62]{rugh}). As continuity of the partials imply differentiability, it follows that~$\Phi$ is continuous and therefore bounded on~${[0,T]^2}$. Consequently, for any~${i,j\in\{1,2,\dots,n\}}$
\[
\sup_{t\in[0,T]}\int_{0}^{t}|\Phi_{ij}(t,\tau)|\,d\tau,\quad \sup_{\tau\in[0,T]}\int_{\tau}^{T}|\Phi_{ij}(t,\tau)|\,dt,
\]
are finite, where~$\Phi_{ij}$ is the entry at the~$i$-th row,~$j$-th column of~$\Phi$. By~\cite[Theorem~75]{vidyasagar}, it follows that the convolution operator is~$\mathcal{L}_p$ stable for all~${p\in[0,\infty]}$.
\begin{rem}
	The bounded integral conditions for~$\mathcal{L}_p$ stability given in~\cite{vidyasagar} are modified here so that the supremums are taken over~$[0,T]$. This is because the continuous state-transition matrix~$\Phi$ can be continuously extended from the compact domain~$[0,T]^2$ to the first quadrant of~$\real^2$~(the system is causal) and ensure a decay fast enough so the conditions hold over an infinite horizon.
\end{rem}


\section{Proof of Proposition~\ref{prop:sufficiency}}
\label{app:sufficiency}
By~\eqref{eq:operator}, the output at pass~${k+1}$ can be written as
\begin{equation*}
		y_{k+1}	=\bar{H}x_{k+1}(0)+\bar{G}_0y_k+\Omega(\varphi(x_{k+1},y_k)),
\end{equation*}
so
\begin{equation}
	y_{k}=\bar{G}_0^ky_0+\sum_{i=1}^{k}\bar{G}_0^{k-i}(\bar{H}x_{i}(0)+\Omega(\varphi(x_{i},y_{i-1})))
\label{eq:perturbsolution}
\end{equation}
for all~${k\in\nnum}$, when the solution exists. Recalling the fact that~${\LInfNorm{\bar{G}_0^k}\leq \bar{M}\bar{\zeta}^k}$ for all~${k\in\nnum}$ for some~${\bar{M}\geq1}$ and~${\bar{\zeta}\in(0,1)}$, from~\eqref{eq:perturbsolution}, it follows that
\begin{multline*}
	\LInfNorm{y_{N}}\leq\bar{M}\bar{\zeta}^N\LInfNorm{y_0}+\max\left\{\LInfNorm{\bar{H}},\LInfNorm{\Omega}\right\}\\
	\times\left(\norm{\mathbf{x}(0)}_{e_1}+\max_{i\in\{1,2,\dots,N\}}\LInfNorm{\varphi(x_{i},y_{i-1})}\right)\sum_{i=1}^{N}\bar{M}\bar{\zeta}^{N-i},
\end{multline*}
therefore
\begin{multline}\label{eq:r1bound}
	\LInfNorm{y_{N}}\leq\underbrace{\bar{M}\bar{\zeta}^N}_{r_1<1}\LInfNorm{y_0}\\
	+\underbrace{\bar{M}\frac{1-\bar{\zeta}^N}{1-\bar{\zeta}} \max\{\LInfNorm{\bar{H}},\LInfNorm{\Omega}\}}_{r_2>0}\\
	\times\left(\norm{\mathbf{x}(0)}_{e_1}+\max_{i\in\{1,2,\dots,N\}}\LInfNorm{\varphi(x_{i},y_{i-1})}\right).
\end{multline}
The rest of the proof will be divided into three steps:


\subsection{Lyapunov Stability}

This part follows the same basic ideas of~\cite[Lemma~3]{cdc2015}. Take any~${\epsilon\in(0,(1-r_1)/r_2)}$, where~$r_1,r_2$ are defined in~\eqref{eq:r1bound}. By Lemmas~\ref{lem:O(x)} and~\ref{lem:finiteL}, there exist~${\delta^*\in(0,\epsilon)}$ and~${\delta_{\rm fh}^*\in(0,\min\{\delta_{\rm fh},\delta^*/L_{\rm fh}\})}$ such that~${\LInfNorm{y_0}+\norm{\mathbf{x}(0)}_{e_1}<\delta_{\rm fh}^*<\delta^*/L_{\rm fh}<\delta^*}$ means
\begin{equation}\label{eq:inbetween}
	\LInfNorm{y_k}\leq L_{\rm fh}(\LInfNorm{y_0}+\norm{\mathbf{x}(0)}_{e_1})<\delta^*<\epsilon
\end{equation}
for all~${k\in\{0,1,\dots,N-1\}}$, which in turn implies
\begin{equation*}
	\begin{aligned}
	\LInfNorm{\varphi(x_{k},y_{k-1})}	&<\epsilon/(L_{\rm fh}+1)(\LInfNorm{y_{k-1}}+\norm{\mathbf{x}(0)}_{e_1})\\
																		&<\epsilon(\LInfNorm{y_0}+\norm{\mathbf{x}(0)}_{e_1})
	\end{aligned}
\end{equation*}
for all~${k\in\{1,2,\dots,N\}}$. Assume~${\LInfNorm{y_0}<\delta_y\leq\delta_{\rm fh}^*/2}$ and~${\norm{\mathbf{x}(0)}_{e_1}<\delta_x\leq r_y\delta_y}$ for arbitrary~$r_y$ satisfying
\[
	r_y\in\left(0,\min\left\{1,\frac{1-r_1-r_2\epsilon}{r_2(1+\epsilon)}\right\}\right).
\]
The interval above is nonempty since~${\epsilon<(1-r_1)/r_2}$, and if~$r_y$ belongs to this interval,~${\delta_x+\delta_y<2\delta_y\leq\delta_{\rm fh}^*}$. It follows by the above arguments and~\eqref{eq:r1bound} that
\begin{multline*}
	\LInfNorm{y_{N}}\leq r_1\LInfNorm{y_{0}}\\
	+r_2(\norm{\mathbf{x}(0)}_{e_1}+\epsilon(\LInfNorm{y_{0}}+\norm{\mathbf{x}(0)}_{e_1}))\\
	\leq\LInfNorm{y_{0}}(r_1 +r_2\epsilon) +\norm{\mathbf{x}(0)}_{e_1}r_2(1+\epsilon),
\end{multline*}
so~${\LInfNorm{y_{N}}\leq\delta_y(r_1 +r_2\epsilon) +\delta_xr_2(1+\epsilon)=r_N\delta_y<\delta_y}$, where~${r_N\triangleq(r_1 +r_2\epsilon) +r_yr_2(1+\epsilon)<1}$. Moreover, by~\eqref{eq:inbetween},~${\LInfNorm{y_k}<\epsilon}$ for all~${k\in\{1,2,\dots,N-1\}}$. By induction,~${\LInfNorm{y_0}<\delta_y}$ and~${\norm{\mathbf{x}(0)}_{e_1}<\delta_x}$ imply~${\LInfNorm{y_k}<\epsilon}$ for all~${k\in\nnum}$, since~${\delta_y<\epsilon}$. Therefore, if~${\LInfNorm{y_0}+\norm{\mathbf{x}(0)}<\delta_1=\min\{\delta_x,\delta_y\}}$, then~${\LInfNorm{y_k}<\epsilon}$ for all~${k\in\nnum}$. As we can find such a~${\delta_1>0}$ for arbitrarily small~${\epsilon>0}$, we conclude that the system is stable.


\subsection{Asymptotic Stability}

From~\eqref{eq:perturbsolution},~$y_{k}	=\bar{y}_k+\sum_{i=1}^{k}\bar{G}_0^{k-i}\Omega(\varphi(x_{i},y_{i-1}))$. Let~${\epsilon=(1-\bar{\zeta})/(2\bar{M}\LInfNorm{\Omega})}$. Since the system is stable, by Lemma~\ref{lem:O(x)} there exists a positive scalar~$\delta_1$ so that~${\LInfNorm{y_0}+\norm{\mathbf{x}(0)_{e_1}}<\delta_2=\delta_1}$ implies
\begin{multline}\label{eq:limsupbound}
	\limsup_{k\to\infty} \LInfNorm{y_k} \\
	\leq \epsilon\bar{M}\LInfNorm{\Omega}\limsup_{k\to\infty} \sum_{i=1}^{k}\bar{\zeta}^{k-i}(\norm{x_{i}(0)}+\LInfNorm{y_{i-1}})\\
	= \epsilon\bar{M}\LInfNorm{\Omega}\limsup_{k\to\infty} \underbrace{\sum_{i=1}^{k}\bar{\zeta}^{k-i}\LInfNorm{y_{i-1}}}_{\bar{S}_k},
\end{multline}
as~${\bar{y}_k\to 0}$, and~${\sum_{i=1}^{k}\bar{\zeta}^{k-i}\norm{x_{i}(0)}\to 0}$ if~${\mathbf{x}(0)\in c_0}$, as we have shown before in Section~\ref{sec:linear}. Now, it is easy to verify that~${\bar{S}_{k+1}=\bar{\zeta}\bar{S}_k+\LInfNorm{y_k}}$, where~$\bar{S}_k$ is defined in~\eqref{eq:limsupbound}. Hence by~\eqref{eq:limsupbound} and Claim~\ref{claim:asymp} we can show
\[
		\limsup_{k\to\infty} \LInfNorm{y_k}\leq\frac{1}{2}\limsup_{k\to\infty} \LInfNorm{y_k},
\]
so~${\limsup_{k\to\infty} \LInfNorm{y_k}=\lim_{k\to\infty} \LInfNorm{y_k}=0}$. Therefore, the system is asymptotically stable.


\subsection{Exponential Stability} 

Let~${\mathbf{x}_{\kappa}(0)\triangleq\{x_{k+1}(0)\}_{k=\kappa}^{\infty}}$ for any~${\kappa\in\nnum}$. As we have proved Lyapunov stability, given~${\epsilon>0}$, by~\eqref{eq:r1bound} and Lemmas~\ref{lem:O(x)} and~\ref{lem:finiteL}, we can find a constant~${\delta_3\in(0,\delta_{\rm fh})}$ such that~${\LInfNorm{y_{0}}+\norm{\mathbf{x}(0)}_{e_\lambda}<\delta_3}$ implies~${\LInfNorm{y_{k}}<\delta_{\rm fh}}$ and
\begin{multline*}
	\LInfNorm{y_{(k+1)N}}\leq r_1\LInfNorm{y_{kN}}\\
	+r_2(\lambda^{kN}\norm{\mathbf{x}(0)}_{e_\lambda}+\epsilon(\LInfNorm{y_{kN}}+\lambda^{kN}\norm{\mathbf{x}(0)}_{e_\lambda}))\\
	\leq(r_1 +r_2\epsilon)\LInfNorm{y_{kN}} +r_2(1+\epsilon)\norm{\mathbf{x}(0)}_{e_\lambda}\lambda^{kN},
\end{multline*}
where we use the~$e_{\lambda}$ norm shift and~$\lambda$ properties, and~$r_1,r_2$ are defined in~\eqref{eq:r1bound}; hence,
\begin{multline*}
	\LInfNorm{y_{kN}}\leq (r_1 +r_2\epsilon)^k\LInfNorm{y_{0}}\\
	+r_2(1+\epsilon)\norm{\mathbf{x}(0)}_{e_\lambda}\sum_{i=1}^k(r_1 +r_2\epsilon)^{k-i}(\lambda^N)^{i-1},
\end{multline*}
for all~${k\in\nnum}$. Now take any
\[
	\epsilon\in\left(\max\left\{0,\frac{1-r_1-2r_2}{3r_2}\right\},\frac{1-r_1}{r_2}\right).
\]
Then,~${r_1+r_2\epsilon<1}$. Letting~${\underline{\lambda}_N\triangleq\max\{r_1+r_2\epsilon,\lambda^N\}}$, as before in the linear case of Section~\ref{subsec:nonzero}, we can find continuous increasing functions
\begin{equation*}
	\begin{aligned}
		K_N(\lambda^N)&\triangleq \max\left\{1,\frac{2r_2(1+\epsilon)}{1-\underline{\lambda}_N}\right\}=\frac{2r_2(1+\epsilon)}{1-\underline{\lambda}_N},\\
		\gamma_N(\lambda^N)&\triangleq\frac{1+\underline{\lambda}_N}{2}\in[\lambda^N,1),
	\end{aligned}
\end{equation*}
by Claim~\ref{claim:exp}, such that~${\LInfNorm{y_{0}}+\norm{\mathbf{x}(0)}_{e_\lambda}<\delta_3}$ implies
\[
	\LInfNorm{y_{kN}}\leq K_N(\lambda^N)\gamma_N(\lambda^N)^k(\LInfNorm{y_{0}}+\norm{\mathbf{x}(0)}_{e_\lambda}),
\]
and since~${\LInfNorm{y_k}\leq\delta_{\rm fh}}$ for all~$k\in\nnum$, by Lemma~\ref{lem:finiteL} and the~$e_{\lambda}$ shift property,
\begin{multline*}
	\LInfNorm{y_{k}}\leq L_{\rm fh}K_N(\lambda^N)\gamma_N(\lambda^N)^{\bar{k}}(\LInfNorm{y_{0}}+\norm{\mathbf{x}(0)}_{e_\lambda})\\
	+L_{\rm fh}(\lambda^N)^{\bar{k}}\norm{\mathbf{x}(0)}_{e_\lambda},
\end{multline*}
for all~${k\in\nnum}$ as~${L_{\rm fh}\geq 1}$, where~$\bar{k}\in\nnum$ satisfies~${k=\bar{k}N+j}$ and~${j\in\{0,1,\dots,N-1\}}$. In turn, this means that
\[
	\LInfNorm{y_{k}}\leq 2L_{\rm fh}K_N(\lambda^N)\gamma_N(\lambda^N)^{\bar{k}}(\LInfNorm{y_{0}}+\norm{\mathbf{x}(0)}_{e_\lambda}),
\]
for all~${k\in\nnum}$. Let~${\gamma(\lambda)\triangleq(\gamma_N(\lambda^N))^{1/N}}$. Then,
\[
	\LInfNorm{y_{k}}\leq 2L_{\rm fh}K_N(\lambda^N)\gamma(\lambda)^{k-j}(\LInfNorm{y_{0}}+\norm{\mathbf{x}(0)}_{e_\lambda}),
\]
hence, as~${\gamma(\lambda)\in(0,1})$ and~${j\leq N-1}$,
\begin{multline}\label{eq:final}
	\LInfNorm{y_{k}}\leq \overbrace{2L_{\rm fh}K_N(\lambda^N)\gamma(\lambda)^{1-N}}^{K(\lambda)}\gamma(\lambda)^k\\
	(\LInfNorm{y_{0}}+\norm{\mathbf{x}(0)}_{e_\lambda}),
\end{multline}
for all~${k\in\nnum}$. Clearly,~$\gamma$ is continuous and increasing as before, while~$K$ defined in~\eqref{eq:final} is continuous. It remains to show that~$K$ is increasing. Since
\begin{equation*}
	\begin{aligned}
	K(\lambda)	&=2L_{\rm fh}\frac{K_N(\lambda^N)}{\gamma(\lambda)^N}\gamma(\lambda)	=2L_{\rm fh}\frac{K_N(\lambda^N)}{\gamma_N(\lambda^N)}\gamma(\lambda)\\
							&=2L_{\rm fh}\frac{2r_2(1+\epsilon)}{1-\underline{\lambda}_N}\frac{2}{1+\underline{\lambda}_N}\gamma(\lambda)\\
							&=8L_{\rm fh}r_2(1+\epsilon)\frac{\gamma(\lambda)}{1-\underline{\lambda}_N^2},
	\end{aligned}
\end{equation*}
it follows that~$K$ is increasing, as~${(1-\underline{\lambda}_N^2)^{-1}}$ is increasing on~${(0,1)}$ as a function of~$\underline{\lambda}_N$.


\begin{ack}
This work was supported by the NSF grant CMMI-1334204, and conducted while the first author was with the Department of Electrical Engineering and Computer Science at the University of Michigan.
\end{ack}


\bibliographystyle{IFAC}
\bibliography{bibs/automatica2017}

\end{document}